\documentclass[11pt]{article}
\usepackage[english]{babel}

\usepackage{a4wide}
\usepackage{parskip}

\usepackage{amsmath}
\usepackage{amsfonts}
\usepackage{latexsym}
\usepackage{graphicx}
\usepackage{xcolor}

\usepackage{bm}

\newcommand{\Nat}{{\mathbb N}}
\newcommand{\Real}{{\mathbb R}}
\newcommand{\Z}{{\mathbb Z}}

\newcommand{\mmod}[1]{\!\!\pmod{#1}}

\def\diag{\mathrm{diag}}

\newcommand{\Cay}{\mathrm{Cay}}
\newcommand{\ord}{\mathrm{ord}}
\newcommand{\diam}{\mathrm{diam}}
\def\vecu{\mbox{\boldmath$u$}}

\def\veca{\mbox{\boldmath$a$}}

\def\vecm{\mbox{\boldmath$m$}}

\newcommand{\HH}{\mathcal{H}}

\newcommand{\LL}{\mathrm{L}}
\newcommand{\sq}[1]{[\![{#1}]\!]}
\newcommand{\sg}[1]{\langle{#1}\rangle}

\newcommand{\cc}{\mathrm{c}}

\def\M{{\bm{M}}}

\def\S{\bm{S}}
\def\U{\bm{U}}
\def\V{\bm{V}}
\def\e{\bm{e}}
\def\a{\bm{a}}
\def\b{\bm{b}}
\def\llambda{\bm{\lambda}}
\def\o{\bm{0}}
\def\x{\bm{x}}
\newcommand{\norm}[1]{\|#1\|_1}

\newtheorem{rem}{Remark}
\newtheorem{lem}{Lemma}

\newtheorem{defi}{Definition}
\newtheorem{teo}{Theorem}

\newtheorem{exa}{Example}
\newtheorem{conj}{Conjecture}

\begin{document}

\title{%
On solid density of\\
Cayley digraphs on finite Abelian groups
      }

\author{
 F. Aguil\'{o} and M. Zaragoz\'{a}\\%
 \\ {\small Departament de Matem\`atiques}
\\ {\small Universitat Polit\`ecnica de Catalunya}
\\ {\small Jordi Girona 1-3 , M\`odul C3, Campus Nord }
\\ {\small 08034 Barcelona. }
\\ {\small\tt \{francesc.aguilo,marisa.zaragoza\}@upc.edu}
}

\date{\today}

\maketitle

\begin{abstract}
Let $\Gamma=\Cay(G,T)$ be a Cayley digraph over a finite Abelian group $G$ with  respect the generating set $T\not\ni0$. $\Gamma$ has order $\ord(\Gamma)=|G|=n$ and degree $\deg(\Gamma)=|T|=d$. Let $k(\Gamma)$ be the diameter of $\Gamma$ and denote $\kappa(d,n)=\min\{k(\Gamma):~\ord(\Gamma)=n,\deg(\Gamma)=d\}$.

We give a closed expression, $\ell(d,n)$, of a tight lower bound of $\kappa(d,n)$ by using the so called {\em solid density} introduced by Fiduccia, Forcade and Zito.

A digraph $\Gamma$ of degree $d$ is called {\em tight} when $k(\Gamma)=\kappa(d,|\Gamma|)=\ell(d,|\Gamma|)$ holds. Recently, the {\em Dilating Method} has been developed to derive a sequence of digraphs of constant solid density. In this work, we use this method to derive a sequence of tight digraphs $\{\Gamma_i\}_{i=1}^{\cc(\Gamma)}$ from a given tight digraph $\Gamma$. Moreover, we find a closed expression of the cardinality $\cc(\Gamma)$ of this sequence. It is perhaps surprising that $\cc(\Gamma)$ depends only on $n$ and $d$ and not on the structure of $\Gamma$.
\end{abstract}



\section{Introduction}\label{sec:intro}

Cayley digraphs $\Gamma=\Cay(G,T)$ over finite Abelian groups $G$ (with generating set $T$) have been used as a model of interconnection networks. Attention has been paid to their diameter $k(\Gamma)$ and its optimization with respect to the order $\ord(\Gamma)=|G|$ and degree $\deg(\Gamma)=|T|$, that is
\begin{equation}
\kappa(d,n)=\min\{k(\Gamma):~\ord(\Gamma)=n,\deg(\Gamma)=d\}.\label{eq:kappa}
\end{equation}
In particular, it is worth studying a finite closed expression of a tight lower bound for the values $\kappa(d,n)$, denoted by $\ell(d,n)$. Usually, the values $\kappa(d,n)$ are obtained by computer search. The more values shared by $\ell(d,n)$ and $\kappa(d,n)$, the better the expression $\ell(d,n)$. As far as we know, $\ell(d,n)$ is only known for $d=1$ and $d=2$, that is $\ell(1,n)=n-1$ (a directed ring with $n$ nodes) and $\ell(2,n)=\lceil\sqrt{3n}\rceil-2$ that was found using a geometrical approach \cite{FYAV:87,WC:74}.

Fiduccia, Forcade and Zito in 1998 \cite{FiFoZi:98} defined the {\em solid diameter} of $\Gamma$ as $D=k+d$, where $d=\deg(\Gamma)$ and $k=k(\Gamma)$. This is the diameter of a {\em minimum distance diagram} $\HH$ (MDD for short) related to $\Gamma$. These diagrams are used to study metric properties of $\Gamma$, mainly the diameter. MDDs and their properties are discussed in the next section. The {\em solid density} of $\Gamma$ is defined to be the density of an MDD related to $\Gamma$, that is $\delta(\Gamma)=\frac{n}{(k+d)^d}$ whith $n=\ord(\Gamma)$. In this work we show that solid density plays a main role in finding $\ell(d,n)$. More precisely, fixed $d$, the {\em global solid density} is defined by $\Delta_d=\sup\{\delta(\Gamma):~\deg(\Gamma)=d\}$. We show (Theorem~\ref{teo:ell}) that $\ell(d,n)=\left\lceil\sqrt[d]{\frac{n}{\Delta_d}}\,\right\rceil-d$.

For a fixed degree $d$, the {\em tightness} $t(d,\Gamma)$ of a digraph $\Gamma$ of degree $d$ and order $n$ is given by $t(d,\Gamma)=k(\Gamma)-\ell(d,n)$. A digraph $\Gamma$ is called {\em tight} when $t(d,\Gamma)=0$. Tightness related results can help in the search for optimal diameter digraphs. Recently, a method \cite{AgFiPe:16c} has been proposed to obtain an infinite family of dense digraphs $\mathcal{F}=\{\Gamma_m\}_{m\geq1}$ generated from an initial dense digraph $\Gamma_1$. The idea of the method is to take an MDD related to $\Gamma_1$, $\HH_1$, and dilate it in a certain way. Then, the family of dilates of $\HH_1$, $\{\HH_m\}_{m\geq1}$, is used to obtain a family of dense digraphs $\mathcal{F}$. The digraphs in $\mathcal{F}$ are called the dilates of $\Gamma_1$. This method is stated in Theorem~\ref{teo:dilM}.

In this work we study the dilates of a given digraph from the point of view of tightness. We see that the diameter of the dilates of $\Gamma_1$ worsens when $t(d,\Gamma_1)>0$ (Theorem~\ref{teo:tight-0}). The case of tight $\Gamma_1$ is fully studied and gives two different cases:
\begin{itemize}
\item[(a)] All the infinite dilates of $\Gamma_1$ are tight. Tight digraphs $\Gamma_1$ with this property are characterized (Theorem~\ref{teo:tight-1}).
\item[(b)] $\Gamma_1$ is tight and only a finite number, $\cc(\Gamma_1)$, of (consecutive) dilates of $\Gamma_1$ are tight. The expression of $\cc(\Gamma_1)$ is found in terms of the order $n=\ord(\Gamma_1)$ and the degree $d=\deg(\Gamma_1)$ (Theorem~\ref{teo:tight-2}).
\end{itemize}

Section~\ref{sec:notation} is devoted to notation and known results. New results for general degree $d$ are given in Section~\ref{sec:dgen}. New complementary results for degrees $d=2$ and $d=3$ are included in Section~\ref{sec:annotd2d3}.

\section{Notation and known results}\label{sec:notation}

Consider a finite Abelian group of order $n$, $G=\Z_{s_1}\oplus\cdots\oplus\Z_{s_d}$ with $n=s_1\cdots s_d$ and $s_1\mid s_2\mid\cdots\mid s_d$, and a {\em generating set} $T=\{g_1,\ldots,g_d\}\subset G$. Sometimes the notation $G=\sg{g_1,\ldots,g_d}$ is used. The {\em Cayley digraph} of $G$ with respect to $T$ is denoted by $\Gamma=\Cay(G,T)$. It has the set of vertices $V(\Gamma)=G$ and the set of arcs $A(\Gamma)=\{g\to g+t:~g\in G, t\in T\}$ and it is strongly connected. The {\em degree} and {\em diameter} of $\Gamma$ are denoted by $\deg(\Gamma)=d$ and   $k(\Gamma)$.

An {\em isomorphism of digraphs} $\psi:D_1(V_1,A_1)\longrightarrow D_2(V_2,A_2)$ is a bijection on the set of vertices $\psi:V_1\rightarrow V_2$ that preserves arcs, i.e. $v_1\to v_2$ belong to $A_1$ iff $\psi(v_1)\to\psi(v_2)$ belongs to $A_2$. An isomorphism between $D_1$ and $D_2$ is denoted by $D_1\cong D_2$.

Consider an integral matrix $\M\in\Z^{d\times d}$ with $n=|\det\M|$, with \textit{Smith normal form} decomposition $\S=\diag(s_1,\ldots,s_d)=\U\M\V$, for unimodular matrices $\U,\V\in\Z^{d\times d}$. Let us denote the Abelian group $G_\M=\Z^d/\M\Z^d$, with the equivalence relation $\a\sim\b$ whenever there is some $\llambda\in\Z^d$ with $\a-\b=\M\llambda$. It is well known that
\begin{equation}
\Cay(G_\M,E_d)\cong\Cay(\Z_{s_1}\oplus\cdots\oplus\Z_{s_d},\{\vecu_1,\ldots,\vecu_d\}),\label{eq:iso}
\end{equation}
where $E_d=\{\e_1,\ldots,\e_d\}$ is the canonical basis of $\Z^d$ and $\vecu_i=\U\e_i$ for $1\leq i\leq d$.

Let be $[r,s)=\{x\in\Real:r\leq x<s\}$. Given $\a=(a_1,\ldots,a_d)\in\Z^d$, we denote the unitary cube $\sq{\a}=\sq{a_1\ldots,a_d}=[a_1,a_1+1]\times\cdots\times[a_d,a_d+1]\subset\Real^d$. The cube $\sq{\a}$ represents the vertex $(a_1,\ldots,a_d)$ in $\Z_{s_1}\oplus\cdots\oplus\Z_{s_d}$. We denote $\sq{\a}\sim\sq{\b}$ whenever $\a\sim\b$ and $\sq{\a}\not\sim\sq{\b}$ otherwise. This equivalence is sometimes denoted as $\a\equiv\b\mmod{\M}$. For a given pair $\a,\b\in\Z^d$, we write $\a\leq\b$ when the inequality $a_i\leq b_i$ holds for each coordinate $1\leq i\leq d$. Let $\Nat$ denote the set of nonnegative integers. Given $\a\in\Nat^d$, consider the set of unitary cubes $\nabla(\a)=\{\sq{\b}:\o\leq\b\leq\a\}$. Given $\x\in\Z^d$, let us consider the $\ell_1$ norm $\norm{\x}=|x_1|+\cdots+|x_d|$.

\begin{defi}[Minimum distance diagram]\label{defi:MDD}
Given a finite Abelian group $G=\sg{g_1,\ldots,g_d}$ of order $n$, consider the map $\phi:\Nat^d\longrightarrow G$ given by $\phi(\a)=a_1g_1+\cdots+a_dg_d$. A {\em minimum distance diagram} (MDD), $\HH$, related to the Cayley digraph $\Gamma=\Cay(G,\{g_1,\ldots,g_d\})$ is a set of $n$ unitary cubes $\HH=\{\sq{\a_0},\ldots,\sq{\a_{n-1}}\}$ such that
\begin{itemize}
\item[$(i)$] $\{\phi(\a):~\sq{\a}\in\HH\}=G$,
\item[$(ii)$] $\sq{\a}\in\HH\Rightarrow\nabla(\a)\subset\HH$.
\item[$(iii)$] $\norm{\a}=\min\{\norm{\x}:~\x\in\phi^{-1}(\phi(\a))\}$ for all $\sq{\a}\in\HH$.
\end{itemize}
\end{defi}

The \textit{diameter of the minimum distance diagram} $\HH$ is defined as $k(\HH)=\max\{\norm{\a}:~\sq{\a}\in\HH\}$. It is well known that for the usual definition of the diameter $k(\Gamma)$ of a Cayley digraph $\Gamma=\Cay(G,\sg{g_1,\ldots,g_d})$, we have $k(\Gamma)=k(\HH)-d$ for every MDD $\HH$ related to $\Gamma$. When $G$ is a cylic group, Definition~\ref{defi:MDD} is equivalent to \cite[Definition~2.1]{SS:09} for multiloop networks. It is also well known that $\HH$ tessellates $\Real^d$ by translation through $d$ independent vectors $S=\{\vecm_1,\ldots,\vecm_d\}$. Let us assume that $\M$ is the $d\times d$ integral matrix of the column vectors in $S$. Then, we have the isomorphisms
\begin{equation*}
\Gamma\cong\Cay(G_{\M},E_d)\cong\Cay(\Z_{s_1}\oplus\cdots\oplus\Z_{s_d},\{\vecu_1,\ldots,\vecu_d\}),
\end{equation*}
where $\{\vecu_1,\ldots,\vecu_d\}$ are the same column vectors of \eqref{eq:iso}. These isomorphisms have already been used in the literature, see for instance \cite{EAF:93,FYAV:87}. 

\begin{defi}[Properness]
Consider a digraph $\Gamma=\Cay(G,T)$ with $G=\Z_{s_1}\oplus\cdots\oplus\Z_{s_d}$. The set $T$ is a {\em proper generating set} when $T=\{\vecu_1,\ldots,\vecu_d\}$, the same set of vectors defined by $\vecu_i=\U\e_i$ in \eqref{eq:iso}.
\end{defi}

Consider a group $G=\Z_{s_1}\oplus\cdots\oplus\Z_{s_d}$ and denote  $tG=\Z_{ts_1}\oplus\cdots\oplus\Z_{ts_d}$. Let us assume that $T$ is a proper generating set of $G$. From the identity $t\S=\U(t\M)\V$, it follows that $T$ is also a (proper) generating set of $tG$.


\begin{defi}[Dilates of MDDs]
Given an MDD $\HH$ related to $\Gamma$, the $m$--dilate of $\HH$, $m\HH$, is defined by $m\HH=\{m\sq{\veca}:~\sq{\veca}\in\HH\}$, where the $m$--dilate of the unitary cube $\sq{\veca}$ is
\[
m\sq{\veca}=\{\sq{m\veca+(\alpha_1,\ldots,\alpha_d)}:~0\leq\alpha_1,\ldots,\alpha_d\leq m-1\}.
\]
\end{defi}

Let us assume that $\HH$ is an MDD related to $\Gamma=\Cay(G,T)$ where $T$ is a proper generating set of $G$. Then, it can be shown that the dilate $t\HH$ is also an MDD related to $t\Gamma=\Cay(tG,T)$. See \cite{AgFiPe:16b,AgFiPe:16c} for more details. Figure~\ref{fig:H} at page \pageref{fig:H} shows an MDD $\HH$ related to the digraph $\Gamma=\Cay(\Z_{16},\{1,4,5\})\cong\Cay(\Z_1\oplus\Z_1\oplus\Z_{16},\{(0,0,1),(0,1,-12),(1,0,-11)\})$ and $2\HH$ related to $2\Gamma\cong\Cay(\Z_2\oplus\Z_2\oplus\Z_{32},\{(0,0,1),(0,1,-12),(1,0,-11)\})$.

\begin{teo}[The Dilating Method \cite{AgFiPe:16b,AgFiPe:16c}]\label{teo:dilM}
Consider a Cayley digraph of degree $d$,\linebreak $\Gamma=\Cay(\Z_{s_1}\oplus\cdots\oplus\Z_{s_d},T)$ where $T=\{g_1,\ldots,g_d\}$ is a proper generating set. Then, for any integer $m\geq1$,
\begin{itemize}
\item[(a)] $\HH$ is an MDD related to $\Gamma$ $\Leftrightarrow$ $m\HH$ is an MDD related to $m\Gamma$.
\item[(b)] $k(m\Gamma)=m(k(\Gamma)+d)-d$.
\end{itemize}
\end{teo}

This method will be referred to as TDM for short. TDM takes advantage of the geometric nature of the minimum distance diagrams. Less is known of these diagrams for degree $d\geq3$. For instance, a generic geometrical description of MDDs is not known for $d\geq3$. One advantage of Theorem~\ref{teo:dilM} is that we can work with the help of MDDs without knowing this generic description. Another advantage is the use of the same generating set for all members of the infinite family of digraphs $\{m\Gamma\}_{m\geq1}$. Indeed, this property comes from the  properness of $T$ as generating set of the initial digraph $\Gamma$.



\section{Main results for general degree}\label{sec:dgen}

First result in this section is a closed tight lower bound $\ell(d,n)$ for the optimal diameter $\kappa(d,n)$ introduced in Section~\ref{sec:intro}. To this end, given a fixed degree $d$, consider the {\em global solid density} $\Delta_d$ defined by
\begin{equation}
\Delta_d=\sup\{\delta(\Gamma):~d(\Gamma)=d\}.\label{eq:Dd}
\end{equation}

\begin{rem}
This notion is well defined. Indeed, from \cite[Theorem~9.1]{DF:04} there is some constant $c$ such that (for $d>1$)
\[
\frac{c}{d(\ln d)^{1+\log_2 e}}\frac{k^d}{d!}+O(k^{d-1})\leq N(d,k)<{d+k\choose d},
\]
where $N(d,k)=\max\{\ord(\Gamma):~\deg(\Gamma)=d,\diam(\Gamma)=k\}$. Then, it follows that
\[
\delta(\Gamma)\leq\frac{N(d,k)}{(k+d)^d}<\frac1{d!k!}<1
\]
and the set $\{\delta(\Gamma):~d(\Gamma)=d\}$ has a supremum $\Delta_d$.
\end{rem}

We consider two types of degrees depending on the fact that $\Delta_d$ is a global maximum or a supremum only.

\begin{defi}
Consider a fixed degree $d\geq1$. The degree $d$ is called {\em closed} when the supremum $\Delta_d$ is a global maximum, otherwise it is called {\em open} degree.
\end{defi}

For a fixed open degree $d$, there is no digraph of degree $d$, $\Gamma$, such that $\delta(\Gamma)=\Delta_d$. In this case it can be assumed the existence of a sequence of digraphs $\{\Gamma_k\}_{k=1}^\infty$ of increasing orders $|\Gamma_k|=n_k$ such that $\lim_{k\to\infty} n_k=\infty$ and $\lim_{k\to\infty}\delta(\Gamma_k)=\Delta_d$. Given a fixed degree $d$, being open or closed, consider a digraph $\Gamma$ with $\deg(\Gamma)=d$ and $\ord(\Gamma)=n$. Since $\delta(\Gamma)\leq\Delta_d$, the inequality
\[
k(\Gamma)\geq\left\lceil\sqrt[d]{\frac{n}{\Delta_d}}\,\right\rceil-d
\]
is always fulfilled. Thus,
\begin{equation}
\kappa(d,n)\geq\left\lceil\sqrt[d]{\frac{n}{\Delta_d}}\,\right\rceil-d\label{eq:kleql}
\end{equation}
holds for each $d$ and $n$.

\begin{teo}\label{teo:ell}
For a given fixed degree $d$ we have
\begin{equation}
\ell(d,n)=\left\lceil\sqrt[d]{\frac{n}{\Delta_d}}\,\right\rceil-d.\label{eq:ell}
\end{equation}
\end{teo}
\noindent{\bf Proof}: Assume $d$ is a closed degree. Assume $\Gamma^*$ is a digraph with $\deg(\Gamma^*)=d$, $\ord(\Gamma^*)=n$ and $\delta(\Gamma^*)=\Delta_d$. Then, $k(\Gamma^*)=\sqrt[d]{\frac{n}{\Delta_d}}-d$ holds. Since $\sqrt[d]{\frac{n}{\Delta_d}}=k(\Gamma^*)+d\in\Nat$, identity $\kappa(d,n)=\left\lceil\sqrt[d]{\frac{n}{\Delta_d}}\,\right\rceil-d$ also holds and the statement expression $\ell(d,n)$ is a tight lower bound for $\kappa$ for this degree.

Assume now that $d$ is an open degree. Now $\delta(\Gamma)<\Delta_d$ for all digraph $\Gamma$ with $\deg(\Gamma)=d$. By definition of the supremum $\Delta_d$, for all $\varepsilon>0$ there is some digraph $\Gamma^*$ with $\deg(\Gamma^*)=d$ such that $0<\Delta_d-\varepsilon<\delta(\Gamma^*)<\Delta_d$. Assume $\ord(\Gamma^*)=n$. From $\Delta_d-\varepsilon<\delta(\Gamma^*)$, it follows that $k(\Gamma^*)<\sqrt[d]{\frac{n}{\Delta_d-\varepsilon}}-d$. On the other hand, inequality $\sqrt[d]{\frac{n}{\Delta_d-\varepsilon}}-d<\sqrt[d]{\frac{n}{\Delta_d}}-d+1$ holds whenever $\varepsilon<\Delta_d-\frac{n}{\left(1+\sqrt[d]{\frac n{\Delta_d}}\right)^d}=\alpha_{d,n}$. Then, choosing $\varepsilon<\min\{\Delta_d,\alpha_{d,n}\}$, we have
\[
\left\lceil\sqrt[d]{\frac{n}{\Delta_d}}\,\right\rceil-d\leq\kappa(d,n)\leq k(\Gamma^*)<\sqrt[d]{\frac{n}{\Delta_d}}-d+1.
\]
Therefore, identity $\kappa(d,n)=\left\lceil\sqrt[d]{\frac{n}{\Delta_d}}\,\right\rceil-d$ also holds and the statement expression $\ell(d,n)$ is a tight lower bound for open degrees also. $\square$

From this last result, it is worth computing the value $\Delta_d$ for each fixed degree $d$. As far as we know, only $\Delta_1$ and $\Delta_2$ are known.

\begin{rem}
Since $\Delta_1=1$ (direct rings of diameter one unit less than the order) and $\Delta_2=\frac13$ (\cite{FiFoZi:98}), expression \eqref{eq:ell} generalizes those known sharp lower bounds for degrees $d=1$ and $d=2$, that is $\ell(1,n)=n-1$ and $\ell(2,n)=\lceil\sqrt{3n}\rceil-2$ (\cite{FYAV:87,WC:74}).
\end{rem}

\begin{rem}\label{rem:Ndk}
Expression \eqref{eq:ell} also gives a tight upper bound $N(d,k)$ for the order a Cayley digraph $\Gamma$ of a finite Abelian group can have
\begin{equation}
|\Gamma|\leq N(d,k)=\left\lfloor\Delta_d(k+d)^d\right\rfloor.\label{eq:Ndk}
\end{equation}
\end{rem}

Assume $d$ is a fixed degree. Given a digraph $\Gamma$, with $\ord(\Gamma)=n$ and optimal diameter $k(\Gamma)=\kappa(d,n)$, can we decide this diameter is good enough? All what we can say is that $\Gamma$ does his best for this particular order value $n$. Expression $\ell(d,n)$ allows us to expand the local goodness for the diameter. This is the idea of {\em tightness}.

\begin{defi}
The {\em tightness} of a digraph $\Gamma$ of degree $\deg(\Gamma)=d$, order $\ord(\Gamma)=n$ and diameter $k(\Gamma)$ is $t(d,\Gamma)=k(\Gamma)-\ell(d,n)$. We say that $\Gamma$ is $t(d,\Gamma)$--tight. Those $0$--tight digraphs are called tight ones.
\end{defi}

For a closed degree $d$, digraphs reaching the maximum density $\Delta_d$ have to be tight. Unfortunately, being tight is not a sufficient condition for a digraph to attain $\Delta_d$. It is worth studying the tightness of dilates generated by TDM with respect to the tightness of the initial digraph.

\begin{teo}\label{teo:tight-0}
Let $\Gamma$ be a non-tight digraph of degree $d$. Let $\mathcal{F}=\{m\Gamma:~m\geq1\}$  be the family of dilates generated by TDM for the initial digraph $\Gamma$. Then, the tightness $t(d,m\Gamma)$ worsens as $m$ grows.
\end{teo}
\noindent{\bf Proof}: Let us assume that $\Gamma$ of order $n$ and degree $d$ is $r$-tight with $r\geq1$. Applying Theorem~\ref{teo:dilM}-(b) to the dilate $\Gamma_m=m\Gamma$, we obtain $k(m\Gamma)=m(k(\Gamma)+d)-d$. Then,
\begin{align*}
k(m\Gamma)&=m\left(\left\lceil\sqrt[d]{\frac{n}{\Delta_d}}\,\right\rceil+r\right)-d=m\left\lceil\sqrt[d]{\frac{n}{\Delta_d}}\,\right\rceil-d+mr\geq\left\lceil m\sqrt[d]{\frac{n}{\Delta_d}}\,\right\rceil-d+mr\\
&=\left\lceil\sqrt[d]{\frac{m^dn}{\Delta_d}}\,\right\rceil-d+mr=\ell(d,|m\Gamma|)+mr.
\end{align*}
Therefore, it follows that $t(d,m\Gamma)\geq mr$. $\square$

As it is stated in Theorem~\ref{teo:tight-0}, the dilating method has not to be used to obtain small diameter digraphs when the initial digraph $\Gamma$ is not tight. Let us study now the behaviour of the dilates when the initial digraph $\Gamma$ is tight. To this end, a characterization of tight dilates is needed.

\begin{lem}[Characterization of tight dilates]\label{lem:charT}
Given a fixed degree $d$, let us assume that $\Gamma$ is a tight digraph with $\deg(\Gamma)=d$ and $\ord(\Gamma)=n$. Then, the dilate $m\Gamma$ is tight iff $\left\lceil m\sqrt[d]{\frac{n}{\Delta_d}}\,\right\rceil=m\left\lceil\sqrt[d]{\frac{n}{\Delta_d}}\,\right\rceil$ holds.
\end{lem}
\noindent{\bf Proof}: Its is a direct consequence of the definition of tightness and Theorem~\ref{teo:dilM}-(b). $\square$

Let us consider the set $C_d$ defined by
\begin{equation}
C_d=\{x\in\Real:~\Delta_dx^d\in\Nat\}.\label{eq:CDelta}
\end{equation}
Let us define $\{x\}=\lceil x\rceil-x$.

\begin{lem}\label{lem:infam}
For a fixed degree $d$, let $\Gamma$ be a digraph with $\deg(\Gamma)=d$ and $\ord(\Gamma)=n$. The identity $\left\lceil m\sqrt[d]{\frac{n}{\Delta_d}}\,\right\rceil=m\left\lceil\sqrt[d]{\frac{n}{\Delta_d}}\,\right\rceil$ holds for all $m\geq1$ iff $n=\Delta_d x^d$ for some $x\in  C_d\cap\Nat$.
\end{lem}
\noindent{\bf Proof}: Fixed $x\in C_d\cap\Nat$, assume $n=\Delta_dx^d$. Thus $\sqrt[d]{\frac{n}{\Delta_d}}=x\in\Nat$ and so, the statement's identity holds for all $m\geq1$.

Assume now that $m\left\lceil\sqrt[d]{\frac{n}{\Delta_d}}\,\right\rceil=\left\lceil m\sqrt[d]{\frac{n}{\Delta_d}}\,\right\rceil$ holds for all $m\geq1$. If $n\neq\Delta_dx^d$ for any $x\in C_d\cap\Nat$, then $\sqrt[d]{\frac{n}{\Delta_d}}=x\notin\Nat$ and $0<\left\{\sqrt[d]{\frac{n}{\Delta_d}}\right\}<1$ hold. So, there is some large enough $m_0\in\Nat$ such that $m_0\left\{\sqrt[d]{\frac{n}{\Delta_d}}\right\}>1$. Therefore, from $m_0\left\lceil\sqrt[d]{\frac{n}{\Delta_d}}\,\right\rceil=m_0\sqrt[d]{\frac{n}{\Delta_d}}+m_0\left\{\sqrt[d]{\frac{n}{\Delta_d}}\right\}$, it follows that $m_0\left\lceil\sqrt[d]{\frac{n}{\Delta_d}}\,\right\rceil>\left\lceil m_0\sqrt[d]{\frac{n}{\Delta_d}}\,\right\rceil$ which makes a contradiction. $\square$

\begin{teo}\label{teo:tight-1}
Let $\Gamma$ be a tight digraph of degree $d$ and order $n=\Delta_dx^d$ for some $x\in C_d\cap\Nat$. Then, all the elements of the family $\mathcal{F}=\{m\Gamma:~m\geq1\}$ generated by TDM are tight.
\end{teo}
\noindent{\bf Proof}:  The statement follows from Lemma~\ref{lem:charT} and Lemma~\ref{lem:infam}. $\square$

Now it is worth studying those dilates that come from a tight digraph that not fulfills Theorem~\ref{teo:tight-1}. Numerical traces point to a finite number of tight dilates. A closed expression of this number, the {\em tightness coefficient}, is found in the following theorem.

\begin{teo}[Tightness Coefficient]\label{teo:tight-2}
Let $\Gamma$ be a tight digraph with $\deg(\Gamma)=d$ and $\ord(\Gamma)=n\neq\Delta_dx^d$ for all $x\in C_d\cap\Nat$. Then, the dilate $m\Gamma$ is tight iff $\displaystyle 1\leq m\leq\cc(d,n)$, where $\cc(d,n)$ is the {\em tightness coefficient} of $\Gamma$ given by
\[
\cc(d,n)=\begin{cases}
         \beta(d,n)-1&\beta(d,n)\in\Nat,\\
\lfloor\beta(d,n)\rfloor&\beta(d,n)\notin\Nat,
         \end{cases}
\quad\text{with}\quad\beta(d,n)=\frac1{\left\{\sqrt[d]{\frac{n}{\Delta_d}}\right\}}.
\]
\end{teo}
\noindent{\bf Proof}: When $x\notin\Nat$, two cases are considered:
\begin{itemize}
\item[(a)] Assume $\beta(d,n)\notin\Nat$. Then, there is a unique $m_0\in\Nat$ such that $0<m_0<\beta(d,n)<m_0+1$. Thus, for all $m\in\Nat$ with $\frac1{m_0}\leq\frac1m$, inequalities $0<\left\{\sqrt[d]{\frac{n}{\Delta_d}}\right\}<\frac1m$ hold. Therefore,
\begin{equation}
0<m\left\lceil\sqrt[d]{\frac{n}{\Delta_d}}\,\right\rceil-m\sqrt[d]{\frac{n}{\Delta_d}}<1\label{eq:ineq}
\end{equation}
holds and so, we also have $\left\lceil m\sqrt[d]{\frac{n}{\Delta_d}}\,\right\rceil=m\left\lceil\sqrt[d]{\frac{n}{\Delta_d}}\,\right\rceil$. Then, by Lemma~\ref{lem:charT}, the digraph $m\Gamma$ is tight for all $1\leq m\leq m_0=\lfloor\beta(d,n)\rfloor$.

\item[(b)] Assume now $\beta(c,n)=m_1\in\Nat$. Then,
\[
1=m_1\left\{\sqrt[d]{\frac{n}{\Delta_d}}\right\}\Leftrightarrow m_1\left\lceil\sqrt[d]{\frac{n}{\Delta_d}}\,\right\rceil-1=m_1\sqrt[d]{\frac{n}{\Delta_d}}
\]
holds. Thus, $\left\lceil m_1\sqrt[d]{\frac{n}{\Delta_d}}\,\right\rceil<m_1\left\lceil\sqrt[d]{\frac{n}{\Delta_d}}\,\right\rceil$ also holds and $m_1\Gamma$ is not a tight digraph (notice that $m_1\geq2$ as a consequence of the tightness of $\Gamma$). Taking $m_0=m_1-1\geq1$ and $1\leq m\leq m_0$, we have $0<m<m_1=\beta(d,n)$. Thus, condition \eqref{eq:ineq} also holds for all these values of $m$. Hence, by Lemma~\ref{lem:charT}, $m\Gamma$ is tight for all $1\leq m\leq m_0=\beta(d,n)-1$.
\end{itemize}

Let us see now that the dilate $m\Gamma$ is not tight for $m>\cc(d,n)=m_0$. We have in both cases (a) and (b) that, if $m\geq m_0+1$, inequalities $m_0<\beta(c,n)\leq m_0+1\leq m$ hold. Then, it follows that $\frac1m\leq\frac1{m_0+1}\leq\left\{\sqrt[d]{\frac{n}{\Delta_d}}\right\}<\frac1{m_0}$ and so $1\leq m\left\lceil\sqrt[d]{\frac{n}{\Delta_d}}\,\right\rceil-m\sqrt[d]{\frac{n}{\Delta_d}}$ holds. Thus, we have $\left\lceil m\sqrt[d]{\frac{n}{\Delta_d}}\,\right\rceil<m\left\lceil\sqrt[d]{\frac{n}{\Delta_d}}\,\right\rceil$ and the digraph $m\Gamma$ is not tight for $m\geq m_0+1$. $\square$

It is a surprising fact that, for a fixed degree $d$, the number of tight dilates of a given tight digraph $\Gamma$  depends only on $n=\ord(\Gamma)$ and not on the structure of $\Gamma$ itself. See examples \ref{exa:dilcoef} and \ref{exa:dilcoefbis} in the next section.

The following result, only valid for closed degrees $d$, characterizes those orders of tight digraphs with solid density attaining the global maximum $\Delta_d$.

\begin{teo}\label{teo:main}
Assume $d$ is a closed degree. Take a digraph $\Gamma$ with degree $d$ and order $n$. Then,
\begin{equation}
\delta(\Gamma)=\Delta_d\Leftrightarrow\Gamma\textrm{ is tight and }n=\Delta_dx^d\textrm{ for }x\in C_d\cap\Nat.
\end{equation}
\end{teo}
\noindent{\bf Proof}: Assume $\delta(\Gamma)=\Delta_d$. Then, $k(\Gamma)=\sqrt[d]{\frac{n}{\Delta_d}}-d$ with  $\sqrt[d]{\frac{n}{\Delta_d}}=x\in\Nat$. Thus, $k(\Gamma)=\ell(d,n)$ holds and $\Gamma$ is tight. Thus, $n=\Delta_dx^d$ for $x\in C_d\cap\Nat$.

Assume now that $\Gamma$ is tight and $n=\Delta_dx^d$ for some $x\in C_d\cap\Nat$. From $k(\Gamma)=\ell(d,n)=\lceil x\rceil-d=x-d$, it follows that $\delta(\Gamma)=\Delta_d$. $\square$

\begin{rem}\label{rem:td}
Fixed a closed degree $d$, consider the rational value $\Delta_d=\frac{s_d}{q_d}$ with $s_d,q_d\in\Nat$ and $\gcd(s_d,q_d)=1$. Assume $q_d=p_1^{\alpha_1}\cdots p_r^{\alpha_r}$ is the prime decomposition of $q_d$ with $\alpha_i=a_id+b_i$ for $a_i,b_i\in\Nat$ with $a_i\geq0$ and $0\leq b_i<d$ for all $i=1\div r$. Defining
\[
\beta_i=\begin{cases}
a_i&\text{if }b_i=0,\\
a_i+1&\text{if }b_i\neq0,
\end{cases}
\]
for $i=1\div r$, the minimum value of $x\in C_d\cap\Nat$ is $x_d=\Pi_{i=1}^rp_i^{\beta_i}$.
\end{rem}

\section{Some annotations for degrees two and three}
\label{sec:annotd2d3}

Clearly, the degree $d=1$ is a simple case. An optimal diameter digraph of order $n$ is a directed ring $\Gamma$ of diameter $k(\Gamma)=n-1$ and density $\delta(\Gamma)=1$. Thus, $\kappa(1,n)=\ell(1,n)=n-1$ and $d=1$ is a closed degree with $\Delta_1=1$. Notice that $\ell(1,n)$ follows the generic expression \eqref{eq:ell} for $d=1$.

The degree $d=2$ is also a closed degree. This result was stated by Forcade and Lamoreaux in 2000 \cite[Section~4]{FL:00} who proved that $\Delta_2=\frac13$ using MDDs. Then, by  Theorem~\ref{teo:ell}, it follows that $\ell(2,n)=\lceil\sqrt{3n}\,\rceil-2$. This result was pointed out by several authors \cite{EAF:93,FYAV:87,R:96,WC:74} who also used MDDs. In this case, minimum distance diagrams are L-shapes (or rectangles). From Remark~\ref{rem:td} and $t_2=3$, we know that the minimum order for a digraph of degree two to attain $\Delta_2$ is $n=\Delta_2t_2^2=3$.  In fact, taking $\Upsilon_2=\Cay(\Z_3,\{2,1\})$, we have $k(\Upsilon_2)=1$ and the maximum density is attained $\delta(\Upsilon_2)=\frac13$.

There are included here some results that complement those of Section~\ref{sec:dgen}. To this end, we remember  some known facts about L-shapes.

\begin{figure}[h]
\centering
\includegraphics[width=0.5\linewidth]{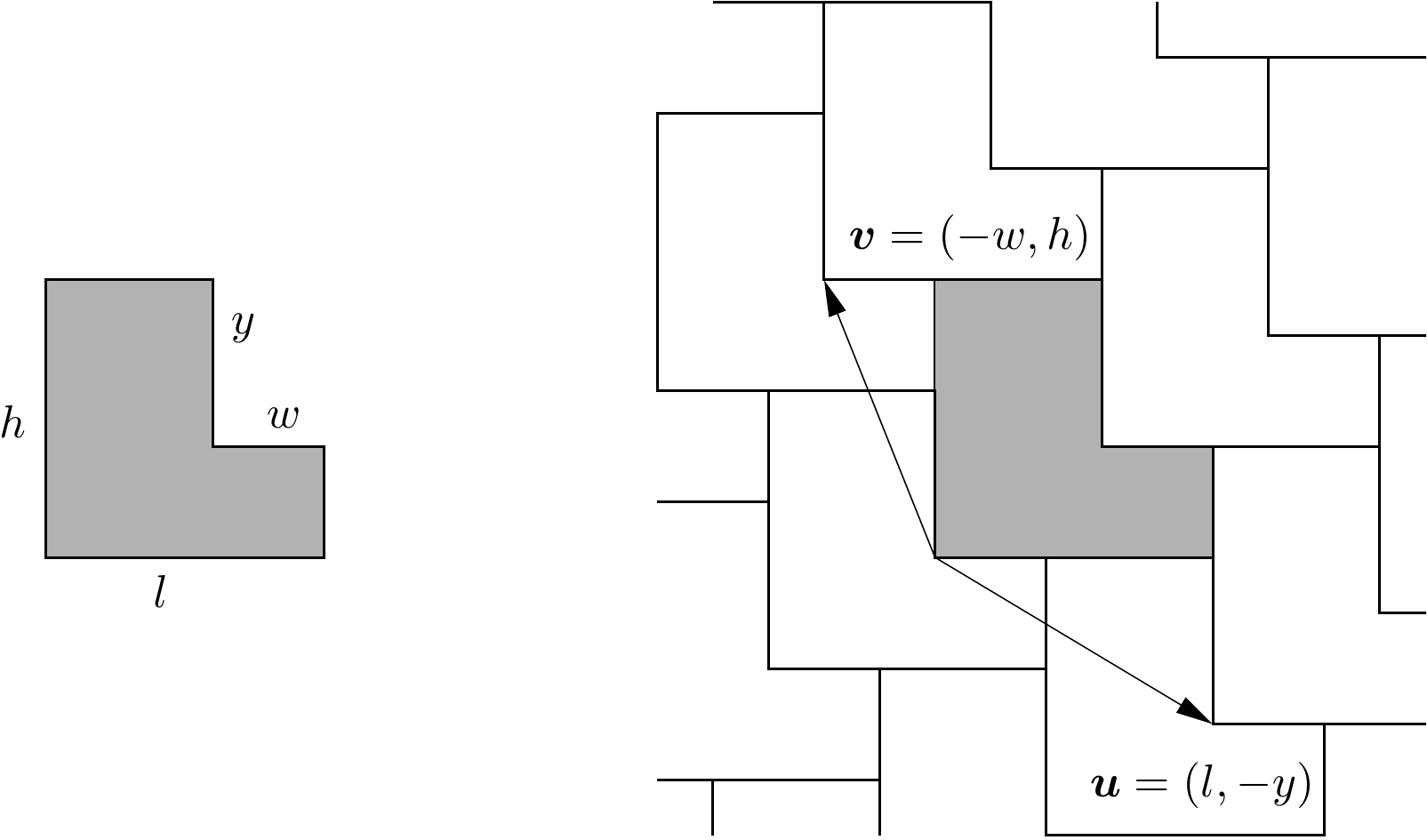}
\caption{Generic L-shape and the related tessellation}
\label{fig:eletes}
\end{figure}

An L-shape $\HH$ related to $\Gamma=\Cay(G,\{a,b\})$ of order $n=s_1s_2$, $s_1\geq1$, with $G=\Z_{s_1}\oplus\Z_{s_2}$ and $s_1\mid s_2$, is denoted by the lengths of its sides, $\HH=\LL(l,h,w,y)$, with $0\leq w<l$ and $0\leq y<h$  (see Figure~\ref{fig:eletes}) and area $n$. Rectangles are particular cases of L-shapes, that is $w=0$ or $y=0$. Let us denote $\gcd(\HH)=\gcd(l,h,w,y)$, $m\HH=\LL(ml,mh,mw,my)$ and $\HH/m=\LL(l/m,h/m,w/m,y/m)$ whenever $m\mid\gcd(\HH)$. Then, $\HH$ is characterized by (see \cite{EAF:93,FYAV:87} for more details)
\begin{align}
n&=lh-wy,\label{eq:area}\\
s_1&=\gcd(\HH),\label{eq:s1}\\
la&=yb\text{ in }G,\label{eq:layb}\\
wa&=hb\text{ in }G,\label{eq:wahb}\\
(l-y)(h-w)&\geq0\text{ and only one factor can vanish},\label{eq:both}
\end{align}
and the solid diameter of $\HH$ (in the sense of Fiduccia, Forcade and Zito \cite{FiFoZi:98}) is given by
\begin{align}
D(\HH)&=l+h-\min\{w,y\}\label{eq:diam}
\end{align}
and the diameter of $\Gamma$  is $k(\Gamma)=D(\HH)-2$. The L-shape $\HH$ tessellates the plane by translation through the vectors $\vecm_1=(l,-y)$ and $\vecm_2=(-w,h)$. Isomorphism \eqref{eq:iso} is given now using the matrix $\M=\left(\begin{array}{rr}l&-w\\-y&h\end{array}\right)$. For instance, the digraph $\Upsilon_2$ has related the L-shape $\HH_2=\LL(2,2,1,1)$ and so, from the Smith normal form decomposition
\[
\diag(1,3)=\left(\begin{array}{rr}0&1\\1&-1\end{array}\right)
\left(\begin{array}{rr}2&-1\\-1&2\end{array}\right)
\left(\begin{array}{rr}1&2\\1&1\end{array}\right),
\]
we have $\Upsilon_2\cong\Cay(\Z_1\oplus\Z_3,\{(0,1),(1,-1)\})$.


As it has been remarked before, the tightness coefficient only depends on the degree and order of the digraph. Below is an example of this fact.

\begin{table}[h]
\small
\centering
\begin{tabular}{|r|r|l|r|r|r|}\hline
$m$&$|m\Gamma|$&$m\Gamma$&L-shape&$k(m\Gamma)$&$\ell(2,|m\Gamma|)$\\\hline\hline
1&72&$\Gamma_1=\Cay(\Z_{72},\{4,11\})$&$\LL(11,8,4,4)$&13&13\\
&&$\Gamma_2=\Cay(\Z_3\oplus\Z_{24},\{(0,1),(-1,3)\})$&$\LL(9,9,3,3)$&&\\\hline
2&288&$2\Gamma_1=\Cay(\Z_2\oplus\Z_{144},\{(-1,4),(-3,11)\})$&$\LL(22,16,8,8)$&28&28\\
&&$2\Gamma_2=\Cay(\Z_6\oplus\Z_{48},\{(0,1),(-1,3)\})$&$\LL(18,18,6,6)$&&\\\hline
3&648&$3\Gamma_1=\Cay(\Z_3\oplus\Z_{216},\{(-1,4),(-3,11)\})$&$\LL(33,24,12,12)$&43&43\\
&&$3\Gamma_2=\Cay(\Z_9\oplus\Z_{72},\{(0,1),(-1,3)\})$&$\LL(27,27,9,9)$&&\\\hline
4&1152&$4\Gamma_1=\Cay(\Z_4\oplus\Z_{288},\{(-1,4),(-3,11)\})$&$\LL(44,32,16,16)$&58&57\\
&&$4\Gamma_2=\Cay(\Z_{12}\oplus\Z_{96},\{(0,1),(-1,3)\})$&$\LL(36,36,12,12)$&&\\\hline
\end{tabular}
\caption{Three dilations of $\Gamma_1$ and $\Gamma_2$ from Example~\ref{exa:dilcoef}}
\label{tab:dilcoef}
\end{table}

\begin{exa}\label{exa:dilcoef}
Let us consider the tight digraphs $\Gamma_1=\Cay(\Z_{72},\{4,11\})\cong\Cay(\Z_1\oplus\Z_{72},\{(-1,4),\allowbreak (-3,11)\})$ and $\Gamma_2=\Cay(\Z_3\oplus\Z_{24},\{(0,1),(-1,3)\})$. Although $\Gamma_1$ has a different structure than $\Gamma_2$, both digraphs are tight and share order. Thus, from $\cc(2,72)=3$, both digraphs have two tight dilations. Table~\ref{tab:dilcoef} has been found by computer.
\end{exa}

The following theorem makes the concept of dilation of a digraph an important tool to be taken into account.

\begin{teo}\label{teo:mainbis}
Let $\Gamma$ be a digraph of degree two that attains the maximum solid density $\Delta_2$. Then, $\Gamma$ is isomorphic to a dilate of $\Upsilon_2$.
\end{teo}
\noindent{\bf Proof}: The degree $d=2$ is closed. Let $\Gamma=\Cay(G,B)$ be a digraph of degree two and solid density $\delta(\Gamma)=\Delta_2$. By Theorem~\ref{teo:main}, $\Gamma$ is tight of order $|\Gamma|=3m^2$. Moreover, $G$ can not be cyclic for $m\geq2$. This fact comes from the characterization of tight Cayley digraphs on finite Abelian groups of degree two by means of L-shapes, included in \cite{EAF:93}. Using this characterization, \cite[Table~2]{EAF:93}, there is only one related tight L-shape $\HH_m=\LL(2m,2m,m,m)=m\HH_2$ of area $3m^2$, where $\HH_2$ is the L-shape related to $\Upsilon_2$. From \eqref{eq:diam}, we have $D(\HH_m)=3m$ for $m\geq1$ (notice that $\HH_m$ fulfills \eqref{eq:area} to \eqref{eq:both}). Then, using the notation \eqref{eq:iso} of Section~\ref{sec:notation}, the related matrix is $\M_m=\left(\begin{array}{rr}2m&-m\\-m&2m\end{array}\right)$, with Smith normal form decomposition
\[
\S_m=\diag(m,3m)=\U_m\M_m\V_m=\left(\begin{array}{rr}0&1\\1&-1\end{array}\right)\M_m\left(\begin{array}{rr}1&2\\1&1\end{array}\right).
\]
By TDM, we obtain the related digraph $\Cay(\Z_m\oplus\Z_{3m},T)=m\Upsilon_2\cong\Gamma$ with the proper generating set $T=\{(0,1),(1,-1)\}$ and diameter $k(\Gamma)=D(\HH_m)-2=3m-2$. $\square$

Contrarily to the case of degree two, less is known about MDDs for degree $d\geq3$. There are no analog to the geometric characterization conditions \eqref{eq:area} to \eqref{eq:diam}. Forcade and Lamoreaux in 2000 \cite{FL:00} proved that the value $\Delta_3'=0.084$ is a local maximum of the solid density for degree $d=3$. We comment here some numerical traces that seem to point that $\Delta'_3$ would be the global solid density for degree three.

As far as we know, no known digraph has solid density larger than $\Delta_3'$. Assuming that $\Delta_3'=\frac{21}{250}$ plays the role of $\Delta_3$, by Theorem~\ref{teo:main} and Remark~\ref{rem:td}, the first tight digraph attaining $\Delta'_3$ would have order $n=\Delta_3'x_3^3=\Delta_3'2^35^3=84$. And this is the case. It is well known that
\[
\Upsilon_3=\Cay(\Z_{84},\{-38,-3,7\})\cong\Cay(\Z_1\oplus\Z_1\oplus\Z_{84},\{(1,10,-38),(0,1,-3),(0,-2,7)\})
\]
is the digraph of smaller order that has solid density $\delta(\Upsilon_3)=\Delta'_3$. Under the previous assumption, there would be an analog to Theorem~\ref{teo:mainbis} for degree three, stating that digraphs attaining $\Delta'_3$ are dilates of $\Upsilon_3$. This infinite family was given in \cite[Proposition~3]{AgFiPe:16c},
\begin{equation}
m\Upsilon_3=\Cay(\Z_m\oplus\Z_m\oplus\Z_{84m},\{(1,10,-38),(0,1,-3),(0,-2,7)\})\label{eq:famd3}
\end{equation}
with diameter $k(m\Upsilon_3)=10m-3$ for $m\geq1$. No other digraph having this density is known. When assuming $\Delta_3=\Delta'_3=0.084$, we write $\ell'(3,n)$ and $\cc'(3,n)$ instead of $\ell(3,n)$ and $\cc(3,n)$, respectively. Then, we also have $\ell'(3,84m^3)=10m-3$ and all these digraphs $m\Upsilon_3$ would be tight.

Remark~\ref{rem:Ndk} and the previous assumption give the following tight upperbound for the order $n$ of Cayley digraphs on Abelian groups of fixed degree three and diameter $k$
\begin{equation}
n\leq N'(3,k)=\left\lfloor\frac{21}{250}(k+3)^3\right\rfloor.\label{eq:Nmaxd3}
\end{equation}
According to this bound, tables 8.1 and 8.2 given by Dougherty and Faber in \cite{DF:04} would give optimal order values for diameters $38$ to $43$ (these orders are marked there as likely optimal values). Expression \eqref{eq:Nmaxd3} agrees with the maximum order when $k\in\{7,8,17,27,37\}$. These values, except $k=8$, correspond to digraphs isomorphic to $m\Upsilon_3$ for $m\in\{1,2,3,4\}$, respectively. The solid density of these digraphs attain  the assumed global maximum $\Delta'_3=0.084$, except the case $k=8$ that has a smaller value $\delta\approx0.83396$.

We can consider an instance to test expression $\cc'(3,n)$ using Theorem~\ref{teo:tight-2}. To this end, we proceed as in Example~\ref{exa:dilcoef}.

\begin{figure}[h]
\centering
\includegraphics[width=0.3\linewidth]{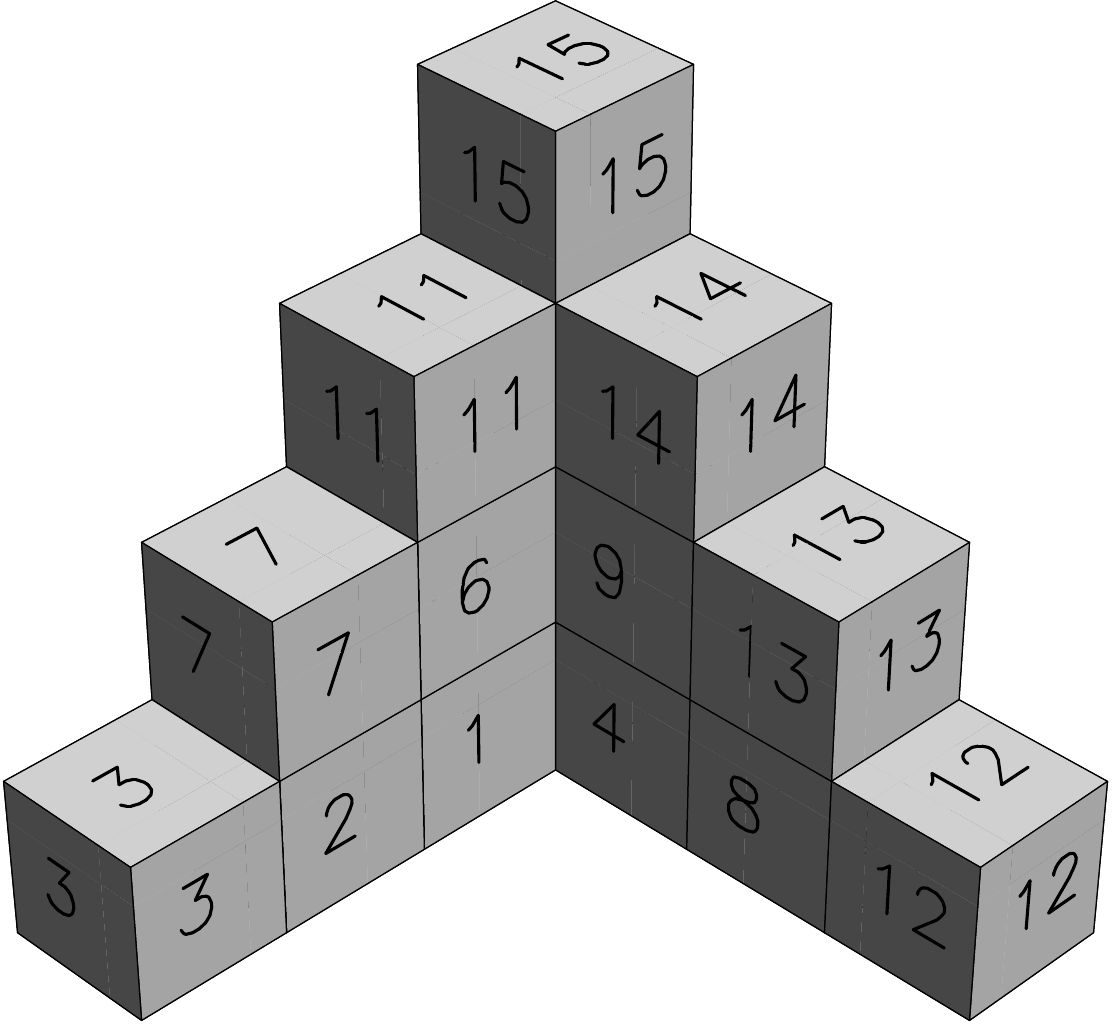}
\hspace{0.01\linewidth}
\includegraphics[width=0.67\linewidth]{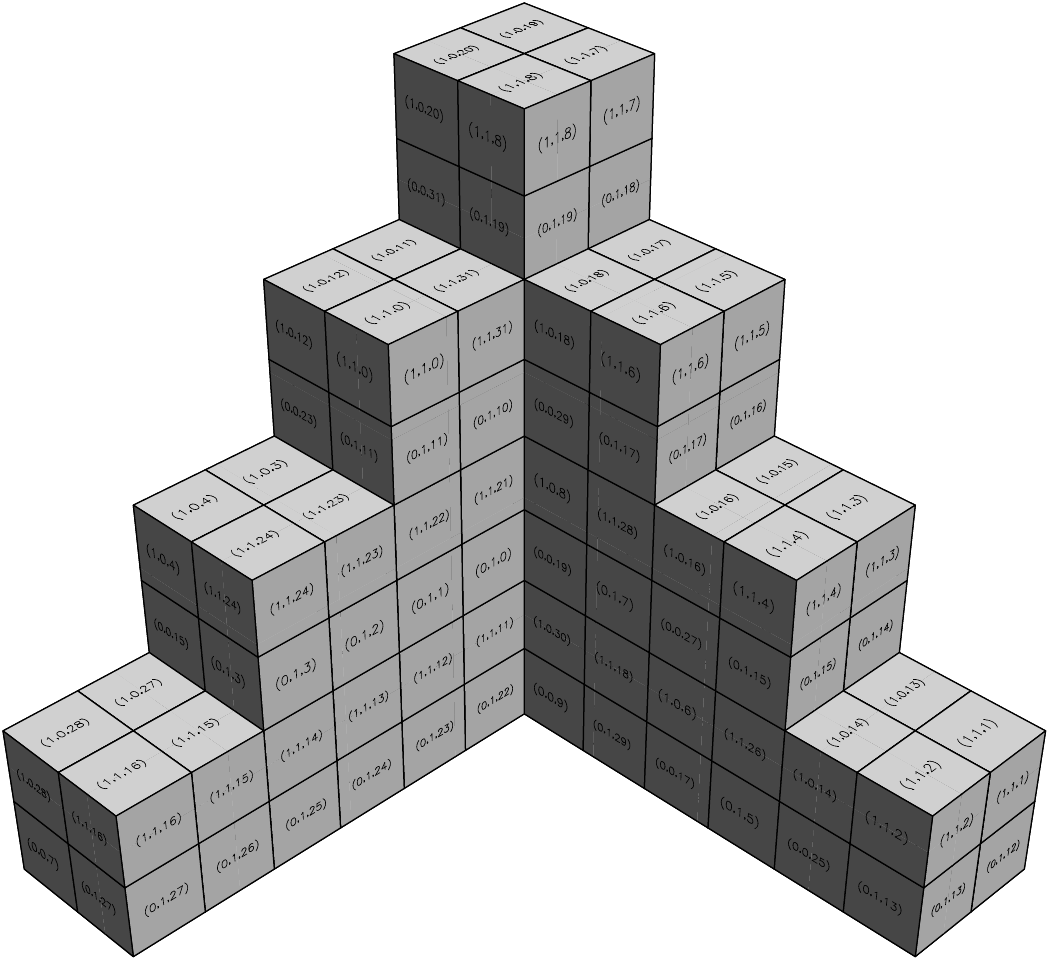}
\caption{MDDs $\HH$ and $2\HH$ related to $\Gamma=\Cay(\Z_{16},\{1,4,5\})$ and $2\Gamma$, respectively}
\label{fig:H}
\end{figure}

\begin{exa}\label{exa:dilcoefbis}
Consider
\[
\Gamma=\Cay(\Z_{16},\{1,4,5\})\cong\Cay(\Z_1\oplus\Z_1\oplus\Z_{16},\{(0,0,1),(0,1,-12),(1,0,-11)\}),
\]
where the isomorphism is given by the Smith normal form decomposition
\[
\diag(1,1,16)=\U\M\V=
\left(\begin{array}{rrr}
0 & 0 & 1 \\
0 & 1 & 0 \\
1 & -12 & -11
\end{array}\right)
\left(\begin{array}{rrr}
-1 & -1 & 0 \\
-1 & 0 & -4 \\
1 & -3 & 0
\end{array}\right)
\left(\begin{array}{rrr}
-8 & -9 & -12 \\
-3 & -3 & -4 \\
2 & 2 & 3
\end{array}\right)
\]
of the matrix $\M$. This is the matrix of column vectors defining the  tessellation of $\Real^3$ by the related MDD $\HH$ (left hand side of Figure~\ref{fig:H}). The MDD $2\HH$ related to $2\Gamma$ is depicted in the right hand side of Figure~\ref{fig:H}. Notice how apparent is the dilation  there, where each cube has been replaced by $2\times2\times2$ cubes. These two MDDs have been obtained by computer search.

\begin{table}[h]
\centering
\begin{tabular}{|r|r|l|r|r|}\hline
$m$&$|m\Gamma|$&$m\Gamma$&$k(m\Gamma)$&$\ell'(3,|m\Gamma|)$\\\hline\hline
1&16&$\Gamma=\Cay(\Z_{16},\{1,4,5\})$&3&3\\
2&128&$2\Gamma=\Cay(\Z_2\oplus\Z_2\oplus\Z_{32},B)$&9&9\\
3&432&$3\Gamma=\Cay(\Z_3\oplus\Z_3\oplus\Z_{48},B)$&15&15\\
4&1024&$4\Gamma=\Cay(\Z_4\oplus\Z_4\oplus\Z_{64},B)$&21&21\\
5&2000&$5\Gamma=\Cay(\Z_5\oplus\Z_5\oplus\Z_{80},B)$&27&26\\\hline
\end{tabular}
\caption{Four dilations of $\Gamma$ with proper set $B=\{(0,0,1),(0,1,-12),(1,0,-11)\}$}
\label{tab:ex3}
\end{table}

The digraph $\Gamma$ of diameter $k(\Gamma)=3$ is tight (according with $\ell'$) and $\cc'(3,|\Gamma|)=4$. Table~\ref{tab:ex3} shows several dilates of $\Gamma$ using the proper generating set $B=\{(0,0,1),(0,1,-12),\allowbreak (1,0,-11)\}$. It can be verified by computer that $\kappa(3,16)=3$, $\kappa(3,128)=9$, $\kappa(3,432)=15$, $\kappa(3,1024)=21$ and $\kappa(3,2000)=27$.
\end{exa}

\begin{figure}[h]
\centering
\includegraphics[width=0.9\linewidth]{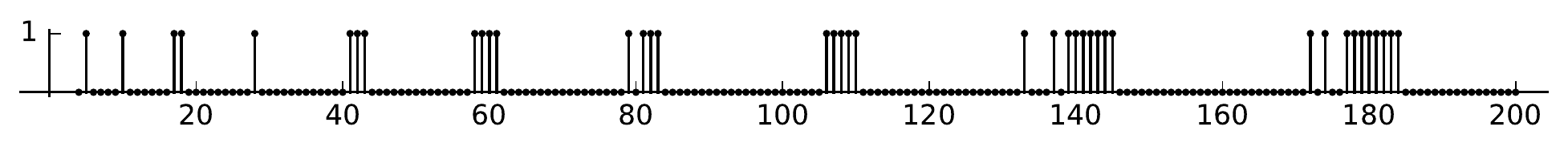}
\caption{$\kappa(3,n)-\ell'(3,n)$}
\label{fig:difer}
\end{figure}

Figure~\ref{fig:difer} shows values of $\kappa(3,n)-\ell'(3,n)$ against $n$ for $4\leq n\leq200$. The maximum value in this range is $1$.

These numerical evidences (and other ones not included here) suggest the following conjecture.

\begin{conj}
The degree $d=3$ is closed and $\Delta_3=\Delta_3'=\frac{21}{250}$.
\end{conj}

Notice that, under this conjecture, the infinite family of digraphs $\Gamma_m=m\Upsilon_3$ of diameter $k(\Gamma_m)=10m-3$, for $m\geq1$, would attain the maximum order $N'(3,10m-3)$ according to \eqref{eq:Nmaxd3}.

\end{document}